\documentclass[11pt, reqno]{amsart}

 \usepackage{amsmath}
 \usepackage{amssymb}
\usepackage{bm}

\DeclareMathAccent{\mathring}{\mathalpha}{operators}{"17}

 \usepackage{color}

\newcommand{\mysection}[1]{\section{#1}
      \setcounter{equation}{0}}

\newtheorem{theorem}{Theorem}[section]
\newtheorem{lemma}[theorem]{Lemma}

\newtheorem{corollary}[theorem]{Corollary} 

\theoremstyle{definition}
\newtheorem{assumption}{Assumption}[section]

\theoremstyle{remark}
\newtheorem{remark}{Remark}[section]

\newtheorem{example}{Example}[section]

\newcommand{\tr}{\text{\rm tr}\,}

 \makeatletter 
 \def\dashint{%  
 \operatorname%
 {\,\,\text{\bf--}\kern-.98em\DOTSI\intop\ilimits@\!\!}}
\def\ninf{\qopname\relax\@empty{inf\phantom{p}\!\!\!}}
 \makeatother

\newcommand\bbeta{\text{\raise-.2ex\hbox{$\bm{\beta}$}}}

\def\sft{{\sf t}}

\newcommand\bR{\mathbb{R}}

\newcommand\cF{\mathcal{F}}

\newcommand\cO{\mathcal{O}}

\newcommand\cN{\mathcal{N}}

\begin{document}

\title[Asimptotic behavior of solutions]
{On the asimptotic behavior
 of solutions of the Cauchy 
problem for parabolic equations
with time periodic coefficients}

\author{R.Z. Khasminskii}
\email{ab3340@wayne.edu}
\address{Wayne State University, Detroit, USA and Institute of the Information Transmission Problems, Moscow, Russia}

\author{N.V. Krylov}
\email{nkrylov@umn.edu}
\address{School of Mathematics, University of Minnesota,
 Minneapolis, MN, 55455, USA}



\begin{abstract}
 We are considering the asimptotic behavior
as $t\to\infty$ of solutions of the Cauchy 
problem for parabolic second order equations
with time periodic coefficients. The problem
is reduced to considering degenerate
time-homogeneous
diffusion processes on the product of a unit circle
and Euclidean space.
  
\end{abstract}

\maketitle

\mysection{Introduction}

Let $\bR^{d}$ be a Euclidean space of points $x=(x^{1},
...,x^{d})$, $d\geq 1$.
We are interested in the asimptotic behavior
as $s\to\infty$ of the solution of the Cauchy
problem
\begin{equation} 
                                      \label{5.14.1}
\partial_{s}u(s,x)=Lu(s,x),\quad 
s>0,x\in\bR^{d},\quad u(0,x)=f(x),
\end{equation}
where $f$ is a given function,
$$
Lu(s,x)=(1/2)a^{ij}(s,x)D_{ij}u(s,x)+b^{i}(s,x)D_{i}u(s,x),
$$
$$
D_{i}=\frac{\partial}{\partial x^{i}},
\quad D_{ij}=D_{i }D_{ j},
\quad \partial_{s}=\frac{\partial}{\partial s}
$$
(the summation convention is enforced throughout the article).

We assume that the coefficients $a=(a^{ij})$ and $b=(b^{i})$ are
$1$-periodic in $t$ and they are
bounded and $a=a^{*}\geq0$. Such seting is quite popular
in many isues, for instance, in parameter
estimation for stochastic processes, the reader is referred to \cite{De_15}, \cite{HK_10}
and the references therein.
Our goal is to find and characterize
$$
\lim_{s\to\infty}u(s,x).
$$
It turns out that, under some assumptions, this limit equals the integral of the product of $f$ and certain
function $p$ which can be found as the value at $s=0$ of the $1$-periodic solution of
$\partial_{s}p+L^{*}p=0$ that is nonnegative and has integral one over $(0,1)\times\bR^{d}$.
 This is, roughly speaking,
the essence of Theorems \ref{theorem 6.14.1}
and \ref{theorem 6.14.2}, see also
Remark \ref{remark 6.22.1}.
Period $1$ is chosen just for convenience,
any other periods can be treated similarly
(also see Example \ref{example 6.29.1}).

The exposition of this result is split into
two parts. The one in Section \ref{section 6.20.2} relies on some qualitative assumptions,
  a conjecture,    
and some  arguments, at the end informal,  which are intuitively
appealing to a probabilist. The other one
in Sections \ref{section 6.20.1}
and \ref{section 6.18.1} is 
completely rigorous albeit
more technical. 

\mysection{Motivation}
                      \label{section 6.20.2}
Here we assume that $a$ and $b$ are
regular enough
to make the constructions in this section 
possible. Suppose that  
for any $(t,x)\in\bR^{d+1}$
there exist  a probability space
and a  $d$-dimensional Wiener process
$w_{t}$ such that
the system
\begin{equation} 
                                      \label{5.14.2}
x _{s}=x  +\int_{0}^{s}\sqrt{a}
( \sft_{r} ,x_{r})\,dw_{r}
+\int_{0}^{s}b ( \sft_{r} ,x_{r}) \,dr,\quad \sft_{s}=
 t-s 
\end{equation}
has a weakly unique solution  $(\sft_{s},x_{s})$  and these
solutions form a time-homogeneous strong  Markov
process in $\bR^{d+1}$. As usual in the
theory of Markov
processes by $E_{t,x}F(\sft_{\cdot},x_{\cdot})$
we mean the expectation of $F(\sft_{\cdot},x_{\cdot})$
where $(\sft_{\cdot},x_{\cdot})$ is the solution 
of \eqref{5.14.2}. The Markov property, in particular,
implies that for any Borel $f(t,x)\geq0$ and any
$s_{1},s_{2}\ge 0$
\begin{equation}
                                      \label{5.14.4}
E_{t,x}E_{\sft_{s_{1}},x_{s_{1}}}f(\sft_{s_{2}},
x_{s_{2}})=
E_{t,x}f(\sft_{s_{1}+s_{2}},x_{s_{1}+s_{2} }).
\end{equation}

In this section we apply to this setting the methods
developed in \cite{Kh} in the case of time-homogeneous diffusion processes. These methods
are applicable when the process possesses
some recurrent properties and, since 
$\sft_{\cdot}$ is not recurrent in any sense
but the coefficients are periodic in time,
it is natural to map $(\sft_{\cdot},x_{\cdot})$
on $\cO\times\bR^{d}$, where $\cO$ is the circle of unit length centered at the origin
on the plane. 

We represent each point in $\cO$
by a number, which is its distance along
$\cO$ to the point $((2\pi)^{-1},0)$ in the   clockwise direction. In this way we define a mapping $y:[0,1)\to\cO$.
We extend this mapping to $\bR$ as $1$-periodic
function. Then even though $y(t)$ is discontinuous
at $0,\pm1,\pm2,...$, if $f(y)$
is, say twice continuously differentiable
on $\cO$, $f(y(t))$ is
twice continuously differentiable
on $\bR$. 

Also introduce the process
  $(y_{s},x_{s}):=(y(\sft_{s}),x_{s})$
on $\cO\times\bR^{d}$. Obviously,
$(y_{s},x_{s})$ is a strong Markov process.
Below the symbol  $E_{y,x}$ stands for the expectation
associated with this process when it starts at
$(y,x)$. The symbol  $P_{y,x}$ is used for
the probability measure
associated with this process when it starts at
$(y,x)$. We also use the symbols $E_{\mu}$ and $P_{\mu}$
when the starting point has distribution $\mu$.
\begin{lemma}
                        \label{lemma 6.14.1}

For $t\in\bR$   and $f\geq0$,
which is Borel on $\cO\times\bR^{d}$, we have
\begin{equation}
                           \label{6.12.1}
E_{y(t),x}f(y_{s},x_{s})=E_{t,x}
f(y(\sft_{s}),x_{s}),
\end{equation}
and if $f=f(y,x)$ is smooth enough
with compact support, then uniformly on
$\bR^{d+1}$
$$
 \frac{1}{s}\big(E_{y(t),x}f(y_{s},x_{s})-f(y(t),x)\big) 
$$
\begin{equation}
                           \label{6.12.4}
\to (1/2)a^{ij}(t,x)D_{ij}
g(t,x)+b^{i}(t,x)D_{i}g(t,x)
-\partial_{t}g(t,x) 
\end{equation}
as $s\downarrow0$, 
where $g(t,x)=f(y(t),x)$.
\end{lemma}

Here \eqref{6.12.1} follows from the definition
of $(y_{s},x_{s})$ and \eqref{6.12.4} is easily
obtained (for instance, under the assumptions 
in Section \ref{section 6.20.1}, see  Lemma \ref{lemma 6.18.1})
by applying It\^o's formula
to $f(y(\sft_{s}),x_{s})=g(\sft_{s},x_{s})$.

Denote $B_{R}=\{x\in\bR^{d}:|x|<R\}$. We impose the following condition (B):

(B1) There exists $R_{0}\in(0,\infty)$ such that the least eigenvalue of $a(t,x)$
is bounded away from zero by a number
independent of $t\in[0,1]$ and $x\in 
B_{R_{0}}$;

(B2) There exists $R_{1}\in(0,R_{0})$
such that, for any $R\geq R_{0}$, we have
$$
\sup_{\cO\times B_{R}}E_{y,x}\tau_{U}<\infty,
$$
where $U=\cO\times \bar B_{R_{1}}$ and 
$\tau_{U}$ is the first time $(y_{s},x_{s})$
hits $U$.

Our first goal is to show that $(y_{s},x_{s})$
has a unique invariant probability
distribution.

Introduce stopping times $0=\tau_{0}<\tau'_{1}
<\tau_{1}...$ so that $\tau'_{n+1}$
is the first time $(y_{s},x_{s})$ hits $\cO\times
\bar B_{R_{1}}$ after $\tau_{n}$,
$\tau_{n+1}$ is the first time $(y_{s},x_{s})$ hits $\cO\times
\partial B_{R_{0}}$ after $\tau'_{n+1}$.
It is not hard to see that,
under condition (B), all $\tau'_{n},\tau_{n}$
are finite (a.s.).

Consider the case that $(y_{0},x_{0})=(y,x)
\in \cO\times \partial B_{R_{0}}$ and
set $(\tilde y_{n},\tilde x_{n})
=(y_{\tau_{n}},x_{\tau_{n}})$. The strong Markov property
of $(y_{s},x_{s})$ implies that
$(\tilde y_{n},\tilde x_{n})$ is a Markov chain on $\cO\times \partial B_{R_{0}}$.

Now comes a subtle point which we could not
justify by referring to the literature known to
us. Under
some natural regularity assumptions on $a$,

(A)
for $(y,x)\in\cO\times\partial B_{R_{1}}$,
the $P_{y,x}$-distribution of $(y_{\tau_{1}},
x_{\tau_{1}})$ has a density bounded away from
zero and infinity by constants independent
of $(y,x)$. 

To discuss this assumption
consider the parabolic equation $\partial_{s}
u=Lu$ in $(0,4]\times B_{R_{0}}$ with zero initial condition and, say regular boundary
condition on $(0,4)\times\partial B_{R_{0}}$.
Usually, $u(s,x)$ is expressed as the integral
over $(0,4)\times\partial B_{R_{0}}$ of
the boundary data agains a (Poisson) kernel $p(s,x,\tau,\xi)$, $(\tau,\xi)\in (0,4)\times\partial B_{R_{0}}$. Having in mind It\^o's formula,
one can easily see that,
for (A) to be true, it suffices to have
an  $\varepsilon>0$ such that
$\varepsilon^{-1}\geq
p(s,x,\tau,\xi)\geq\varepsilon$ whenever
$s\in(3,4)$,
$|x|=R_{1}$, $\tau\in (1,2)$, and $|\xi|=R_{0}$.

{\em Conjecture\/}. Property (A) holds is Assumption \ref{assumption 6.19.1} is  
satisfied.

As in Section 4.4 of \cite{Kh},
in light of (A) and (B) and the result
on page 197 of \cite{Do} (see there (5.6)),
we arrive at the following.

\begin{lemma}
                       \label{lemma 6.19.1}
The Markov chain $(\tilde y_{n},\tilde x_{n})$,
$n=1,2,...$, has a unique stationary 
probability distribution $\tilde \mu$.
Furthermore, for a constant $\kappa\in(0,1)$
and all $n\geq 1$ and Borel 
$\Gamma\subset \cO\times\partial B_{R_{0}}$
\begin{equation}
                             \label{6.19.2}
|\tilde P^{n}(y,x,\Gamma)-
\tilde \mu(\Gamma)|\leq \kappa^{n-1}
\end{equation}
Finally, $\tilde \mu$ has a density
bounded away from zero and infinity.
\end{lemma}

As in Section 4.4 of \cite{Kh},
Lemma \ref{lemma 6.19.1} leads to
existence of invariant probability
distribution for $(y_{t},x_{t})$.
\begin{theorem}
                   \label{theorem 6.19.1}
For Borel $\Gamma\subset \cO\times\bR^{d}$
define
$$
\nu(\Gamma)=\int_{\cO\times\partial B_{R_{0}}}
E_{y,x}\int_{0}^{\tau_{1}}I_{\Gamma}
(y_{s},x_{s})\,ds\,\tilde \mu(dydx).
$$
Then $\bar\nu:=\nu(\cO\times\bR^{d})<\infty$
and $\mu:=\nu/\bar \nu$ is an
invariant probability
distribution for $(y_{t},x_{t})$.

\end{theorem}

To prove this theorem it suffices to repeat
a short and straightforward proof
of Theorem  4.1 of \cite{Kh}.

To prove that the stationary probability
distribution is unique we need
the law of large numbers.

\begin{theorem}
                   \label{theorem 6.19.2}
Let $f$ be a Borel function on
$\cO\times\bR^{d}$ having finite integral
with respect to $\tilde\mu$ from Theorem 
\ref{theorem 6.19.1}. Then for any
$(y,x)\in \cO\times\bR^{d}$ 
\begin{equation}
                             \label{6.19.3}
\frac{1}{T}\int_{0}^{T}f(y_{s},x_{s})
\,ds\to\int_{\cO\times\bR^{d}}f(y',x')\,\mu(dy'dx')
\end{equation}
as $T\to\infty$ $P_{y,x}$-(a.s.).

\end{theorem}

Proof. We follow the   proof
of Theorem 4.2 of \cite{Kh} assuming for
simplicity that $f$ is bounded. First we prove
that  \eqref{6.19.3} holds $P_{\tilde\mu}$-(a.s.).
Introduce
$$
\eta_{n}=\int_{\tau_{n}}^{\tau_{n+1}}f(
y_{s},x_{s})\,ds.
$$
Note that (by definition)
$$
E_{\tilde\mu}\eta_{n}=\int_{\cO\times\bR^{d}}f(y,x)\,\nu(dydx).
$$
Furthermore,
under the measure $P_{\tilde\mu}$ the sequence
$\eta_{n}$ is obviously stationary and,
therefore, there exists $\xi$ such that 
$P_{\tilde\mu}$-(a.s.)  and in $L_{1}$
\begin{equation}
                           \label{6.20.1}
\frac{1}{n}\sum_{i=1}^{n}\eta_{i}\to \xi
\end{equation}
as $n\to\infty$. It follows that for any integer $m$
$$
\lim_{n\to\infty}\frac{1}{n-m}\sum_{i=m}^{n}E_{\tilde\mu} \{\eta_{i}
\mid \cN^{y_{\cdot},x_{\cdot}}_{\tau_{m}}
 \}=E_{\tilde\mu}\{\xi
\mid \cN^{y_{\cdot},x_{\cdot}}_{\tau_{m}}\},
$$
where $\cN^{y_{\cdot},x_{\cdot}}_{t}$
is the $\sigma$-field generated by
$(y_{s},x_{s})$, $s\leq t$.
By the strong Markov property and the stationarity of $(y_{\tau_{n}},x_{\tau_{n}})$
for $i\geq m$
$$
E_{\tilde\mu} \{\eta_{i}
\mid \cN^{y_{\cdot},x_{\cdot}}_{\tau_{m}}
 \}=
E_{\tilde\mu} \{E\{\eta_{i}
\mid \cN^{y_{\cdot},x_{\cdot}}_{\tau_{i}}\}
\mid \cN^{y_{\cdot},x_{\cdot}}_{\tau_{m}}
 \}
$$
$$
=E_{\tilde\mu} \{\eta(y_{\tau_{i}},x_{\tau_{i}})
\mid \cN^{y_{\cdot},x_{\cdot}}_{\tau_{m}}
 \}=\zeta_{i-m}(y_{\tau_{m}},x_{\tau_{m}}),
$$
where $\eta$ and $\zeta_{i-m}$ are
some Borel functions. Hence
$$
\lim_{n\to\infty}\frac{1}{n-m}\sum_{i=m}^{n}\zeta_{i-m}(y_{\tau_{m}},x_{\tau_{m}})=E_{\tilde\mu}\{\xi
\mid \cN^{y_{\cdot},x_{\cdot}}_{\tau_{m}}\},
$$
which is equivalently rewritten as
$$
\lim_{n\to\infty}\frac{1}{n }\sum_{i=0}^{n}\zeta_{i}(y_{\tau_{m}},x_{\tau_{m}})=E_{\tilde\mu}\{\xi
\mid \cN^{y_{\cdot},x_{\cdot}}_{\tau_{m}}\}.
$$
By taking into account that the distribution of
$(y_{\tau_{m}},x_{\tau_{m}})$ is $\tilde\mu$,
we conclude that the limit of $(1/n)
\sum_{i=0}^{n}\zeta_{i}$ exists $\tilde\mu$-(a.e.) and if we take a Borel $\kappa(y,x)$
equal to this limit $\tilde\mu$-(a.e.),
then
$$
E_{\tilde\mu}\{\xi
\mid \cN^{y_{\cdot},x_{\cdot}}_{\tau_{m}}\}=
\kappa(y_{\tau_{m}},x_{\tau_{m}}),
$$
\begin{equation}
                                \label{6.19.5}
\xi=\lim_{m\to\infty}
E_{\tilde\mu}\{\xi
\mid \cN^{y_{\cdot},x_{\cdot}}_{\tau_{m}}\}
=\lim_{m\to\infty}\kappa(y_{\tau_{m}},x_{\tau_{m}}).
\end{equation}

In light of Lemma \ref{lemma 6.19.1}
this yields that
for any integer $r\geq1$
$$
E_{\tilde\mu}\{\xi
\mid \cN^{y_{\cdot},x_{\cdot}}_{\tau_{r}}\}
=\lim_{m\to\infty}E_{\tilde\mu}\{\kappa(y_{\tau_{m}},x_{\tau_{m}})
\mid \cN^{y_{\cdot},x_{\cdot}}_{\tau_{r}}\}
$$
$$
=\int_{\cO\times\partial B_{R_{0}}}
\kappa(y,x)\,\tilde\mu(dydx),
$$
which together with \eqref{6.19.5} imply that
$\xi=\text{const}$. This constant
is easily found from \eqref{6.20.1} to be equal to $E_{\tilde\mu}\eta_{1}$ and therefore,
as $n\to\infty$,
$$
\frac{1}{n}\int_{0}^{\tau_{n+1}}f(y_{s},
x_{s})\,ds\to \int_{\cO\times\bR^{d}}f(y ,x )\,\nu(dy dx )
$$
$P_{\tilde\mu}$-(a.s.). By taking here $f\equiv1$
we find that $\tau_{n+1}/n\to \bar\nu$ and this 
and the boundedness of $f$ readily  imply that \eqref{6.19.3} holds
$P_{\tilde\mu}$-(a.s.), that is, it holds
$P_{y,x}$-(a.s.) for $\tilde\mu$-almost any 
$(y,x)$. Actually, the boundedness of $f$
is not needed, as it is shown by a simple argument in the proof of Theorem 4.2 of \cite{Kh}.

For any other $(y,x)\in\cO\times\bR^{d}$, the $P_{y,x}$-distribution of $(y_{\tau_{1}},x_{\tau_{1}})$
on $\cO\times\partial B_{R_{0}}$ is absolutely continuous
 because of (A)
and the strong Markov property. Hence, it is also absolutely continuous with respect to
$\tilde\mu$. Now the fact that \eqref{6.19.3}
 holds under $P_{y,x}$ follows from
$$
\lim_{T\to\infty}
\frac{1}{T}\int_{0}^{T}f(y_{s},x_{s})
\,ds=\lim_{T\to\infty}
\frac{1}{T}\int_{\tau_{1}}^{\tau_{1}+T}f(y_{s},x_{s})
\,ds.
$$
The theorem is proved.

After that the fact that the stationary
probability distribution for $(y_{s},x_{s})$
is unique is derived as in Corollary 4.4 of \cite{Kh}:
From \eqref{6.19.3} for Borel $\Gamma$
and invariant probability distribution $\tilde
\mu$ of $(y_{s},x_{s})$ we have 
$$
\frac{1}{T}\int_{0}^{T}P(y,x,s,\Gamma)\,ds=
E_{y,x}\frac{1}{T}\int_{0}^{T}I_{\Gamma}(y_{s},x_{s})
\,ds\to\mu(\Gamma),
$$
$$
\tilde\mu(\Gamma)=
\int_{\cO\times\bR^{d}}
\Big(\frac{1}{T}\int_{0}^{T}P(y,x,s,\Gamma)\,ds\Big)
\,\tilde\mu(dydx)\to\mu(\Gamma).
$$

\begin{remark}
                      \label{remark 6.22.2}
With the above $\mu$ on $\cO\times\bR^{d}$ we naturally
associate a measure $\hat\mu$ on
$[0,1)\times\bR^{d}$ by the formula
$$
\hat\mu([s,t)\times \Gamma)
=\mu([y(s),y(t))\times\Gamma),
$$
where $[y(s),y(t))$ is the arc between $y(s)$
and $y(t)$. Then for any Borel nonnegative
$f$ on $\cO\times\bR^{d}$ we have
$$
\int_{\cO\times\bR^{d}}f\,\mu(dydx)
=\int_{[0,1)}\int_{\bR^{d}}f(y(t),x)\,
\hat\mu(dtdx).
$$

Since $\mu$ is the invariant probability
distribution of $(y_{s},x_{s})$,   the
integral of the left-hand side of
\eqref{6.12.4} against $\hat\mu$ is zero
by the above, hence,
\begin{equation}
                               \label{6.22.3}
\int_{[0,1)}\int_{\bR^{d}}\Big[(1/2)a^{ij}(t,x)D_{ij}
g(t,x)+b^{i}(t,x)D_{i}g(t,x)
-\partial_{t}g(t,x)\Big]\,\hat\mu(dtdx)=0
\end{equation}
for any $1$-periodic in time smooth $g$.
It is usual in such situations to call
$\hat\mu$ a  {\em generalized solution\/} of 
\begin{equation}
                               \label{6.22.4}
\partial_{t}\hat\mu+(1/2)D_{ij}(a^{ij}\hat\mu)
-D_{i}(b^{i}\hat\mu)=0 
\end{equation}
in $[0,1)\times\bR^{d}$.

\end{remark}
\begin{remark}
                      \label{remark 6.21.1}
For $f=f(y)$ in \eqref{6.19.3} we find that
$$
\int_{\cO\times \bR^{d}}f(y)\,\mu(dydx)
=\int_{0}^{1}f(y(t))\,dt,
$$
which implies that, if $\cO\times\bR^{d}$-valued
random element $(\eta,\xi)$ has distribution
$\mu$, then $\eta$ is uniformly distributed
on $\cO$. Let $\pi_{y}(\Gamma)$, $y\in\cO$, 
$\Gamma$ Borel in $\bR^d$, be the regular
conditional distribution of $\xi$ given
$\eta=y$. Then for all Borel $f\geq0$
$$
\int_{\cO\times \bR^{d}}f(y,x)\,\mu(dydx)
=\int_{0}^{1}\Big(\int_{\bR^{d}}f(y(t),x)\,
\pi_{y(t)}(dx)\Big)\,dt,
$$
and $\hat\mu(dtdx)=\hat \pi_{t}(dx)dt$,
where $\hat \pi_{t}(dx)=\pi_{y(t)}(dx)$.

Also observe that $y_{1}=y_{0}$ no matter
which initial distribution on 
$\cO\times \bR^{d}$ we take. Therefore,
the invariance of $\mu$ implies that
for $f(y,x)=g(y)h(x)$ with Borel $g,h\geq0$ we have
$$
\int_{\cO\times \bR^{d}}g(y)
E_{y,x}h( x_{1})\,\mu(dydx)=
\int_{\cO\times \bR^{d}}
E_{y,x}f(y_{1},x_{1})\,\mu(dydx)
$$
$$
=\int_{\cO\times \bR^{d}}g(y)h(x)\,\mu(dydx).
$$
The arbitrariness of $g$ implies that 
on $(0,1)$ (a.e.)
$$
\int_{\bR^{d}}E_{t,x}h(x_{1})\,\hat\pi_{t}(dx)
=\int_{\bR^{d}} h(x )\,\hat\pi_{t}(dx),
$$
that is, $\hat\pi_{t}$ is an ivariant
probability distribution for the kernel
$P_{t,x}(x_{1}\in\Gamma)$. Quite often
in such situations
$$
E_{t,x}h(x_{n})\to \int_{\bR^{d}} h(y )\,\hat\pi_{t}(dy)
$$
as $n\to\infty$. Observe that as a rule
$E_{t,x}h(x_{n})$ equals $u(n,x)$, where 
$u(s,x)$
is the solution of the Cauchy problem
$\partial_{s}u=Lu$, $s\geq t$, with initial value
$u(t,y)=h(y)$ and this is a motivation
of the way we pursue the main goal of
the article in the next sections.
\end{remark}

\mysection{Proving that the limit
of solutions exists}
                     \label{section 6.20.1}

We  recall  the main setting and introduce
new assumptions.
\begin{assumption}
                     \label{assumption 6.19.1}
For a fixed $\alpha\in(0,1)$,   $a,b$ are defined in $\bR\times\bR^{d}$,
are  
 $1$-periodic in the 
time variable, are $\alpha$ H\"older continuous in $x$
with the H\"older constant independent of $s$, 
are $\alpha/2$ H\"older continuous in $s$
with the H\"older constant independent of $x$,
and
$a$ is symmetric with eigenvalues
between $\delta$ and $\delta^{-1}$
for a constant $\delta\in(0,1]$.
\end{assumption}

The regularity assumptions can be weakened
but our goal is not to get the most general result but 
rather outline our approach to the problem
and its solution. It is well known that
even under much weaker regularity assumptions
there exists strong Markov process
$(\sft_{s},x_{s})$ the trajectory of
which satisfy \eqref{5.14.2} if $(\sft_{s},x_{s})$ starts at $(s,x)$, where 
$w_{t}$ is a $d$-dimensional Wiener process
(see, for instance,
\cite{SV_79} or \cite{Kr}).

\begin{assumption}
                     \label{assumption 6.19.10}
For an $R_{0}\in(0,\infty)$,
\begin{equation}
                                   \label{5.8.1}
(b(s,x),x/|x|)\leq -\gamma/|x|,\quad  |x|\geq R_{0},
\end{equation}
where the constant $\gamma$ is such that  $2\gamma>\delta^{-2}d$. 
\end{assumption}

\begin{lemma}
                               \label{lemma 5.13.1}
 For any $t\in\bR$, $R\in(0,\infty)$, $s>0$, 
there exists $\varepsilon
=\varepsilon(d,\delta,R,s)>0$ such that for any Borel
$\Gamma\subset \bR^{d}$  
 we have 
$$
\inf_{|x|\leq R}P_{t,x}(x_{s}\in
\Gamma)\geq  e^{-1/\varepsilon}\int_{\Gamma}e^{-
|y|^{2}/\varepsilon}\,dy.
$$

\end{lemma}

Proof. As follows from the theory of parabolic equations
\cite{Fr_64}, under our conditions on $a,b$
for any smooth bounded $f\geq0$
and $T\in\bR$, there exists a smooth
solution of $\partial_{t}u(t,x)=Lu(t,x)$
for $t>T$ with boundary condition $u(T,x)=f(x)$. Applying It\^o's formula to $u(\sft_{s},x_{s})$
yields that $u(t,x):=E_{t,x}f(x_{t-T})$.
On the other hand, there is an integral representation of $u$ by means of the fundamental
solution which implies that  for $t>T$
$$
E_{t,x}f(x_{t-T})\geq \int_{\bR^{d}}f(y)p_{t-T}(y-x)\,dy,
$$
where  $p_{t-T}(y)$ is a gaussian-like function. The arbitrariness of $f$ yields the desired result.
The lemma is proved.

In the terminology from \cite{MT_96}
(see there Proposition 4.2.1 and Section 5.4.3)
Lemma \ref{lemma 5.13.1} implies that  the 
time-homogeneous
Markov chain $x_{n},n=0,1,..$ 
(the spacial coordinate of $(\sft_{n},x_{n})$)
is irreducible
relative to Lebesgue measure, strongly aperiodic, and
all balls $B_{R}:=\{x:|x|<R\}$ are 
small and petite.

\begin{lemma}
                               \label{lemma 5.13.2}
Set $V(x)=|x|^{2}$.
Then there exists  constants $\varepsilon>0$
and $R_{1}\geq R_{0}$
such that for $|x|\geq R_{1}$ and any $t$ 
 we have 
$$
E_{t,x}V( x_{1})\leq V( x)-\varepsilon.
$$ 

\end{lemma}

Proof. We have with $\lambda=x/|x|$
$$
LV(t,x)= \tr a(t,x) 
+ 2(  b(t,x),x)\leq  \delta^{-2}d
+2(  b(t,x),x)\leq \delta^{-2}d-2\gamma=:-\kappa.
$$
By It\^o's formula
$$
E_{t,x}V( x_{1})=V(x)+E_{t,x}\int_{0}^{1}
LV(t-s,x_{s})\,ds
\leq V(x)-\kappa
$$
$$
+E_{t,x}\int_{0}^{1}
\big|LV(t-s,x_{s})+ \kappa\big|
I_{|x_{s}|\leq R_{0}}\,ds.
$$
Since the coefficients are bounded
$$
\sup_{t}E_{t,x}\int_{0}^{1}
\big|LV(t-s,x_{s})+ \kappa\big|
I_{|x_{s}|\leq R_{0}}\,ds\to 0
$$
as $|x|\to \infty$. This yields the desired result.
The lemma is proved.

This lemma implies that condition  (iii)
of Theorem 14.0.1 of \cite{MT_96} is fulfilled
for any $t$ and the transition kernel
$P_{t,x}(x_{1}\in \Gamma)$. According to this theorem
for any $t$ there is a unique invariant distribution
$\pi_{t}$ for this kernel and for any $x$
$$
\sup_{f\in\cF}|E_{t,x}f(x_{n})-(f,\pi_{t})|\to0
$$
as $n\to \infty$, where
$$
(f,\pi_{t})=\int_{\bR^{d}}
f(y)\,\pi_{t}(dy)
$$
and $\cF$ is the set of 
continuous $f$ on $\bR^{d}$
satisfying $|f|\leq1$.
\begin{remark}
                              \label{remark 5.13.1}
Obviously $\pi_{t}$ is $1$-periodic function of $t$   because
for any $(t,x)\in\bR^{d}$, $s\geq0$, $n=0,1,...$, and Borel
bounded $f$ we have
$$
E_{t,x}f(x_{s})=E_{t+n,x}f(x_{s}).
$$
\end{remark}
This equality also implies that
$$
\sup_{f\in \cF}|E_{t+n,x}f(x_{n})-(f,\pi_{t})|\to0 
$$
as $n\to \infty$. 
By It\^o's formula,   $E_{t+n,x}f(x_{n})=u_{t}(n,x)$,
where $u_{t}$ is the  
  solution  of the Cauchy
problem
\begin{equation}
                                  \label{5.13.1}
\partial_{s}u_{t}(s,x)=Lu(s,x),\quad s>t,\quad u_{t}(t,x)=f(x),
\end{equation}
given by means of the fundamental solution
(see, for instance, \cite{Fr_64}).

It follows from the   estimates for
the fundamental solution (\cite{Fr_64}) that, for any $t\in[0,1]$ and $R\in(0,\infty)$, the set of solutions of
\eqref{5.13.1} when $f\in\cF$,  
$s\geq2$, $|x|\leq R$ is equicontinuous as functions
of $(s,x)$. Hence, for any $t\in[0,1]$,
$$
u_{t}(n,x)\to (f,\pi_{t})
$$
as $n\to\infty$
uniformly on any compact set uniformly
with respect to $f\in\cF$. Next, observe that, for
$\tau\in[0,1]$,
$$
|u_{t}(n+\tau,x)-(f,\pi_{t})|
\leq E_{ \tau,x}|u_{t}
(n,x_{\tau})-(f,\pi_{t})|
$$
$$
\leq \sup_{|y|\leq R}|u_{t}
(n,y)-(f,\pi_{t})|+2P_{ \tau,x}(|x_{\tau}|>R).
$$
For any $x$ the last probability goes to zero
uniformly with respect to $\tau\in[0,1]$
since $a$ and $b$ are bounded. 
Thus we arrive at the following.
\begin{theorem}
                        \label{theorem 6.14.2}

For any
$t\in[0,1]$,
$$
u_{t}(s,x)\to (f,\pi_{t})
$$
as $s\to\infty$
uniformly on any compact set uniformly
with respect to $f\in\cF$. 
\end{theorem}

In the next section we need a few properties of $\pi_{t}$.

\begin{lemma}
                              \label{lemma 5.13.20}
If $f\in\cF$, then for any $s\geq0$
\begin{equation}
                            \label{6.14.6}
\int_{\bR^{d}}E_{t,x}f(x_{s})\,\pi_{t}(dx)=\int_{\bR^{d}}
f(x)\,\pi_{t-s}(dx).
\end{equation}
\end{lemma}

Indeed, the left-hand side is the limit of
$$
E_{t+n,y}E_{t,x_{n}}f(x_{s})=
E_{t-s+n,y}f(x_{n}).
$$

\begin{corollary}
                    \label{corollary 6.14.1}
For any bounded Borel $f$ the function $(f,\pi_{t})$
is continuous.

\end{corollary}

Indeed, the left-hand side of \eqref{6.14.6}
is continuous with respect to $s$ if
$f$ is continuous. The case of general $f$
is taken care of by approximation and the following
remark.

\begin{remark}
                       \label{remark 6.4.3}
For any Borel $f\geq 0$ and $t\in \bR$ we have
$$
\int_{\bR^{d}}E_{t,x}f(x_{1})\,\pi_{t}(dx)
=\int_{\bR^{d}} f(x )\,\pi_{t}(dx).
$$
Here, as we know from the theory
of parabolic equations (\cite{Fr_64}),
$$
E_{t,x}f(x_{1})\leq N\int_{\bR^{d}}f(y)
\exp(-|x-y|^{2}/N)\,dy,
$$
where the constant $N$ is independent of $f$, $t$. It follows that $\pi_{t}$
has a {\em bounded\/} density  
$\rho(t,x)\leq N$, where $N$ is independent 
of $t$. The continuity of $(f,\pi_{t})$
just for bounded continuous $f$ implies that
we may choose $\rho(t,x)$ to be a Borel function of
$(t,x)$.

\end{remark}

\mysection{How to find  $\pi_{t}$ analytically?}
                       \label{section 6.18.1}

We need some notation. Recall that $\alpha
\in(0,1)$ is a fixed number. Then
$C^{\alpha/2,\alpha}(\bR^{d+1})$ is defined
as the space of functions $f(t,x)$
that are $\alpha$ H\"older continuous in $x$
with the H\"older constant independent of $s$
and 
are $\alpha/2$ H\"older continuous in $t$
with the H\"older constant independent of $x$.
The space $C^{1+\alpha/2,2+\alpha}(\bR^{d+1})$
consists of functions $f$ such that
$\partial_{t}f, D^{2}f,Df ,f
\in C^{\alpha/2,\alpha}(\bR^{d+1})$.
 
Next, recall that $\cO$ is the 
circle of unit length centered at the origin on the   plane.   In light of Remark \ref{remark 6.4.3},  
 the measure
\begin{equation}
                           \label{6.14.3}
\Pi(\Gamma):=\int_{0}^{1}\int_{\bR^{d}}I_{\Gamma}(y(t),x)
\,\pi_{t}(dx)\,dt  =
\int_{0}^{1}\int_{\bR^{d}}I_{\Gamma}(y(t),x)
\,\rho(t,x)\,dxdt  
\end{equation}
is well defined on Borel subsets $\Gamma$ of 
$\cO\times\bR^{d}$. 
We introduce 
$C^{\alpha/2,\alpha}(\cO\times\bR^{d})$ and
$C^{1+\alpha/2,2+\alpha}(\cO\times\bR^{d})$
similarly to how it   was done above.
\begin{lemma}
                         \label{lemma 6.18.1}
Assertions \eqref{6.12.4} of Lemma \ref{lemma 6.14.1}
holds true if $f\in C^{1+\alpha,2+\alpha}
(\cO\times\bR^{d})$.
\end{lemma}

This is easily proved as it is outlined after
Lemma \ref{lemma 6.14.1} by using the
assumed regularity of $a$, $b$, and the
derivatives of $f$ (and, hence, of $g$).

 Next, for $(y,x)\in\cO\times\bR^{d}$
and Borel bounded $f$ on $\cO\times\bR^{d}$ define
$$
Rf(y,x)=E_{y,x}\int_{0}^{\infty}
e^{-s}f(y_{s} ,x_{s})\,ds.
$$

\begin{lemma}
                                \label{lemma 5.14.1}
The function $Rf(y,x)$ is continuous
on $\cO\times\bR^{d}$ for any
Borel bounded $f$.
\end{lemma}
Proof.
Since for $y=y(t)$, owing to
\eqref{6.12.1}, we have
$$
Rf(y,x)
=E_{t,x}\int_{0}^{\infty}
e^{-s}f(y(t-s),x_{s})\,ds=:\hat Rf(t,x),
$$
it suffices to prove that $\hat Rf(t,x)$
is continuous on $ \bR^{d+1}$ for any Borel bounded~$f$.

The parabolic Aleksandrov estimate implies that
$$
\hat R |f|(t,x)\leq N_{1}\Big(\int_{\bR^{d}}\int_{0}^{\infty}
|e^{-s}f(y(t-s),x)|^{d+1}\,dsdx\Big)^{1/(d+1)}
$$
$$
\leq N_{2}\Big(\int_{\bR^{d}}\int_{0}^{1}
|f(y(t-s),x)|^{d+1}\,dsdx\Big)^{1/(d+1)}
$$
\begin{equation}
                                  \label{5.14.5}
=N_{2}\Big(\int_{\bR^{d}}\int_{0}^{1}
|f(y(s),x)|^{d+1}\,dsdx\Big)^{1/(d+1)},
\end{equation}
where the constants $N_{i}$ are independent of $f$ and
$(t,x)$. Here it is appropriate to observe that, as is easy to see, for any $R\in(0,\infty)$,
$$
\lim_{\rho\to\infty}\sup_{|x|\leq R,t\in[0,1]}
\int_{0}^{\infty}e^{-s} P_{t,x}(|x_{s}|\geq\rho)\,ds=0
$$
and therefore it suffices to consider
only the $f$'s that vanish
as $|x|\geq \rho$ for some $\rho$. 
 
Such an $f$ can be approximated in the 
$L_{d+1}$ norm on $\cO\times\bR^{d}$
by uniformly bounded  $f_{k}$, each of which
is  Lipschitz
continuous and   estimate \eqref{5.14.5}
  shows that in that case we have $\hat R  f _{k}
\to \hat R  f$ uniformly. Therefore we
may assume that $f$ is Lipschitz continuous
and has compact support.
Then observe that
$$
\hat R f(t,x)=\check R f(t,t,x),
$$
where
$$
\check R f(\tau,t,x):
=E_{t,x}\int_{0}^{\infty}
e^{-s}f(y(\tau-s),x_{s})\,ds.
$$
Since $f$ is Lipschitz continuous,
$\check R f(\tau,t,x)$ is continuous in $\tau$
uniformly with respect to $(t,x)$. The fact that
it is continuous (smooth) in $(t,x)$ for any $\tau$
follows  from the fact that, as a function of $(t,x)$
the function $\check R f(\tau,t,x)$ is a unique
 $C^{1+\alpha/2,2+\alpha}(\bR^{d+1})$-solution of
$$
\partial_{t}v(t,x)=Lv(t,x)-v(t,x)
+f(y(\tau-t),x),\quad (t,x)\in\bR^{d+1}.
$$
The lemma is proved.

\begin{theorem}
                            \label{theorem 5.13.1}
The measure $\Pi$ is a unique invariant measure
of the process $(y_{s},x_{s}) $.
\end{theorem}

Proof. Let $ f(y,x)$ be Borel bounded
on $\cO\times\bR^{d}$. 
According to \eqref{6.12.1} and Lemma \ref{lemma 5.13.20}
$$
\int_{\cO\times\bR^{d}}E_{y,x}
f(y_{s},x_{s})\,\Pi(dy,dx)
$$
$$
=\int_{0}^{1}\int_{\bR^{d}}
E_{t,x}
f(y(t-s),x_{s})\,\pi_{t}(dx)\,dt
$$
$$
=\int_{0}^{1}\int_{\bR^{d}}
f(y(t-s),x)\,\pi_{t-s}(dx)\,dt
$$
$$
=\int_{0}^{1}\int_{\bR^{d}}
f(y(t),x)\,\pi_{t}(dx)\,dt=
\int_{\cO\times\bR^{d}} 
f(y ,x )\,\Pi(dy,dx).
$$

This proves that $\Pi$ is indeed an invariant
measure. Furthermore,
$$
\int_{\cO\times\bR^{d}}Rf(y,x)
\,\Pi(dy,dx)=\int_{0}^{\infty}
e^{-s}\int_{\cO\times\bR^{d}}E_{y,x}f(y_{s} ,x_{s})\,\Pi(dy,dx)\,ds 
$$
\begin{equation}
                               \label{6.13.1}
=\int_{0}^{\infty}
e^{-s}\int_{\cO\times\bR^{d}} 
f(y ,x )\,\Pi(dy,dx)\,ds
=\int_{\cO\times\bR^{d}} 
f(y ,x )\,\Pi(dy,dx),
\end{equation}
which proves that $\Pi$ is an invariant
measure for the kernel $R$.

 Then the uniqueness of
 $\Pi$ now follows from
Lemma \ref{lemma 5.14.1} due to
Corollary 2.7 of \cite{Ha_16}.
The theorem is proved.

\begin{remark}
                          \label{remark 6.14.1}
The integral over $(0,1)\times\bR^{d}$
of the left-hand side of
\eqref{6.12.4} with respect to 
$\pi_{t}(dx)dt=\rho(t,x)dxdt$ by definition is $1/s$ times
$$
\int_{\cO\times\bR^{d}}E_{y,x}
f(y_{s},x_{s})\,\Pi(dy,dx)-
\int_{\cO\times\bR^{d}} 
f(y,x)\,\Pi(dy,dx),
$$
which is zero. Therefore,
\begin{equation}
                               \label{6.13.3}
\int_{0}^{1}\int_{\bR^{d}}\Big[(1/2)a^{ij}(t,x)D_{ij}
g(t,x)+b^{i}(t,x)D_{i}g(t,x)
-\partial_{t}g(t,x)\Big]\rho(t,x)\,dxdt=0
\end{equation}
for any $1$-{\em periodic\/} in time $g\in 
C^{1+\alpha/2,2+\alpha}(\bR^{d+1})$.
It is usual in such situations to call
$\rho$ a  {\em generalized solution\/} of 
\begin{equation}
                               \label{6.13.2}
\partial_{t}\rho+(1/2)D_{ij}(a^{ij}\rho)
-D_{i}(b^{i}\rho)=0 
\end{equation}
in $(0,1)\times\bR^{d}$. Obviously,
one can change arbitrarily $\rho$ on a set
of measure zero in $(0,1)\times\bR^{d}$
and still preserve \eqref{6.13.3}. Additional properties that
restrict the choice of generalized solutions of
\eqref{6.13.2} follow
from Corollary \ref{corollary 6.14.1}:

(i) For any $t\in(0,1)$, $\rho(t,x)\geq0$
(a.e.) and
$$
\int_{0}^{1}\int_{\bR^{d}}\rho(t,x)\,dx
 dt 
=1,
$$

(ii) For any bounded continuous $f$
$$
(f,\rho(t,\cdot)):= \int_{\bR^{d}}f(x)\rho(t,x)\,dx
$$
 is
continuous on $(0,1)$.

\end{remark}

We have proved the first assertion of the following.

\begin{theorem} 
                        \label{theorem 6.14.1}
(a) For any $t\in\bR$, $\pi_{t}$ has a density
$\rho(t,x)$ which is a generalized
solution of \eqref{6.13.2} in $(0,1)\times\bR^{d}$,
that possesses the properties (i) and (ii).

(b) Suppose that we found a  
generalized
solution $\bar \rho(t,x)$ of
\eqref{6.13.2} in $(0,1)\times\bR^{d}$,
that possesses the properties (i) and (ii).
Then there exists a unique
(up to (a.e.)) bounded Borel $\bar\rho(0,x)$
 such that
for $\bar\rho(1,x):=
\bar\rho(0,x)$ and any bounded Borel $f$

(1)
the function
$$
 \int_{\bR^{d}}f(x)\bar\rho(t,x)\,dx
$$
is continuous on $[0,1]$,

(2) it
 coincides with  
$(f,\pi_{t}) $ on $[0,1]$.
\end{theorem}

Proof. Let $g\in C^{\alpha/2,
\alpha}(\bR^{d+1})$ be $1$-periodic in $t$,
$\lambda>0$, and
$$
G_{\lambda}(t,x)=E_{t,x}
\int_{0}^{\infty}e^{-\lambda s}g(\sft_{s},x_{s})
\,ds.
$$
Then, as it follows from \cite{Fr_64},  $G$ is of class
$C^{1+\alpha/2,2+
\alpha}(\bR^{d+1})$ and is a unique
solution of
$$
(1/2)a^{ij}(t,x)D_{ij}
G_{\lambda}(t,x)+b^{i}(t,x)D_{i}G_{\lambda}(t,x)
-\partial_{t}G_{\lambda}(t,x)
$$
$$
=\lambda G_{\lambda}(t,x)-g(t,x)
$$
in $\bR^{d+1}$. Due to    periodicity
of $g$ and uniqueness (or obviously), $G_{\lambda}$ is $1$-periodic in $t$.
Then since $\bar \rho(t,x)$ is a generalized
solution  of
\eqref{6.13.2}, \eqref{6.13.3} holds
with $G$ in place of $g$, so that
$$
\int_{0}^{1}\int_{\bR^{d}}[
\lambda G_{\lambda}(t,x)-g(t,x)]\bar\rho(t,x)
\,dxdt=0.
$$
This is rewritten as
$$
\lambda\int_{0}^{\infty}e^{-\lambda s}
\Big(\int_{0}^{1}\int_{\bR^{d}}E_{t,x}g(\sft_{s},x_{s})
\bar\rho(t,x)\,dxdt\Big)ds
$$
\begin{equation}
                               \label{6.13.5}
=
\int_{0}^{1}\int_{\bR^{d}} g(t,x) \bar\rho(t,x)
\,dxdt.
\end{equation}
Since $E_{t,x}g(\sft_{s},x_{s})$ is continuous
in $s$ uniformly with respect to $(t,x)$
and the integral of $\bar\rho$ is finite,
the integrand of the first integral in
the left-hand side of \eqref{6.13.5}
is continuous. Then, since \eqref{6.13.5}
holds for any $\lambda>0$ we have for any $s\geq0$
$$
\int_{0}^{1}\int_{\bR^{d}}E_{t,x}g(\sft_{s},x_{s})
\bar\rho(t,x)\,dxdt=
\int_{0}^{1}\int_{\bR^{d}} g(t,x) \bar\rho(t,x)
\,dxdt,
$$
which means that the probability measure
$$
\bar\Pi(\Gamma):=\int_{0}^{1}\int_{0}^{\infty}
I_{\Gamma}(y(t),x)\bar\rho(t,x)
\,dxdt
$$
is an invariant measure for $(y_{s},x_{s})$.
Its uniqueness implies that $\bar \rho=\rho$
almost everywhere on $(0,1)\times\bR^{d}$.
Then for any bounded continuous $f$,
$(f,\pi_{t})=(f,\bar \rho(t,\cdot))$
for almost all $t\in(0,1)$. Since both parts
are continuous in $t$ the equality holds for
all $t\in(0,1)$, but then the equality
holds for all bounded Borel $f$
and $t\in(0,1)$. Since the limits
of $(f,\pi_{t})$ as $t\downarrow0$ and
$t\uparrow 1$ exist and coincide and
$\rho$ and, hence, $\bar\rho$ are bounded,
there exist bounded $\bar\rho(0,x)$ and
$\bar\rho(1,x)$ for which the assertions
of the theorem hold true. 
The theorem is proved.

\begin{remark}
                     \label{remark 6.22.1}
Fix $f\in\cF$ and recall that $u_{t}(s,x)$ 
is introduced as the solution  of the Cauchy
problem \eqref{5.13.1} with initial data
$u_{t}(t,x)=f(x)$. By Theorems \ref{theorem 6.14.2} and \ref{theorem 6.14.1},
for any $t\in(0,1)$ and
any
  generalized solution $\bar\rho$
of  
\eqref{6.13.2} in $(0,1)\times\bR^{d}$,
that possesses the properties (i) and (ii)
we have
$$
u_{t}(s,x)\to \int_{\bR^{d}}f(z)\bar\rho(t,z)\,dz
$$
as $s\to\infty$
uniformly on any compact set.

\end{remark}

\begin{example}
                          \label{example 6.29.1}
Let $\beta=\beta(s)$ be a smooth function
$l$-periodic function ($l>0$)
on $\bR$ which is strictly positive and such that
$1+\beta'\geq1/2$. Let $d=1$,
 introduce $b(s,x)=-x(1+\beta')/(2\beta) $, $a=2$,
and consider the Cauchy problem
$$
\partial_{s}u=D^{2}u+bDu,\quad s>0,\quad u(0,x)=f(x),
$$
where $f$ is a smooth bounded function on $\bR$
and $D=\partial/(\partial x), D^{2}=DD$.
This setting is different from the above one because
$b$ is unbounded. Still we believe that
Theorems \ref{theorem 6.14.1} and \ref{theorem 6.14.2}
are applicable here
after we reduce the period $l$ to one.
To do that set $u(s,x)=v(s/l,x)$. Then $v$ satisfies
$$
\partial_{s}v(s,x)=lD^{2}v(s,x)+lb(ls,x)D v(s,x)
$$
and $b(ls,x)$ is $1$-periodic in $s$. Equation \eqref{6.13.2} takes the form
$$
\partial_{s}\rho(s,x)+lD^{2}\rho(s,x)-
D(lb(ls,x) \rho(s,x))=0,
$$
which turns out to have a classical and, hence,
generalized solution given by
$$
\rho(s,x)=(4\pi\beta(ls))^{-1/2}\exp\big(-x^{2}/(4\beta(ls))\big).
$$
This function obviously possesses the properties
(i) and (ii), and therefore
$$
\lim_{s\to\infty}u(s,x)=
\lim_{s\to\infty}v(s,x)=(2\pi\beta(0))^{-1/2}\int_{-\infty}^{\infty}f(x)e^{-x^{2}/(2\beta(0))}\,dx.
$$
\end{example}
\medskip

{\bf Acknowledgments}. The authors are sincerely
grateful to Yu. Kutoyants for fruitful discussions
and comments.
  
\medskip

{\bf Conflict of interest}  The authors declare that they have no conflict of interest.

\end{document}